\documentclass[10pt]{article}
\usepackage[dvips]{graphicx}
\usepackage[tbtags]{amsmath}
\usepackage{amssymb,color}

\definecolor{c20}{rgb}{0.,0.7,0.}
\definecolor{c30}{rgb}{0.,0.,1.}
\definecolor{c40}{rgb}{1,0.1,0.7}
\definecolor{c50}{rgb}{1,0,0}

\catcode`\@=11
\@addtoreset{equation}{section}
\catcode`\@=12

\def\R{\mathbb{R}}

\allowdisplaybreaks[4]
\newcommand{\COM}[1]{}

\newcommand{\norm}[1]{\| #1 \|}

\renewcommand{\Box}{\fbox{ \hspace{-1.2mm}}}
\newcommand{\Pb}{{\bf P}}
\newcommand{\pk}[1]{\mbox{\Pb}\{#1\}}
\newcommand{\E}[1]{\mbox{\rm${\bf{E}}$}\{#1\}}
\newcommand{\Eg}[1]{\mbox{\rm${\bf{E}}$}\Bigl\{#1 \Bigr\}}
\def\nj{{\bf 1}}

\newcommand{\sprod}[2]{ \langle #1, #2  \rangle}

\newcommand{\abs}[1]{\lvert #1 \rvert}

\newtheorem{theorem}{Proposition}[section]

\newtheorem{corollary}[theorem]{Corollary}
\newtheorem{lemma}[theorem]{Lemma}

 \def\IF{\infty}

\def\kHya{{\cal H}^0_1}
\def\kHy{{\cal H}_2^0}

\newcommand{\BQN}{\begin{eqnarray}}
\newcommand{\EQN}{\end{eqnarray}}
\newcommand{\BQNY}{\begin{eqnarray*}}
\newcommand{\EQNY}{\end{eqnarray*}}
\newcommand{\BS}{\begin{sat}}
\newcommand{\ES}{\end{sat}}
\newcommand{\BL}{\begin{lemma}}
\newcommand{\EL}{\end{lemma}}
\newcommand{\BT}{\begin{theo}}
\newcommand{\ET}{\end{theo}}
\newcommand{\BK}{\begin{corollary}}
\newcommand{\EK}{\end{corollary}}
\newcommand{\pb}[1]{\mbox{\rm$\vk{P}$}\Bigl \{#1 \Bigr \}}
\newcommand{\vk}[1]{\boldsymbol{#1}}
\def\uhh{\underline{h}}
\def\uhhV{\underline{h}^*}
\def\lhh{\overline{h}}
\def\CON{V}
\def\PCON{\widetilde{\CON}}
\def\POLPRJ{\PCON_{p,h}}
\def\PRJ{\CON_{p,h}}

\def\CONB{W}
\def\PCONB{\widetilde{\CONB}}
\def\POLPRJB{\PCONB_{p,h}}
\def\PRJB{\CONB_{p,h}}

\def\njTAU{\nj_{u}(B_0(s,t))}

\def\deltahH{h- \uhh}
\def\we{\uhh_\epsilon}
\def\BVHA{BV_H([0,1]^2)}

\def\BVVA{BV_V([0,1]^2)}

\begin{document}

\begin{center}
{\Large \bf Boundary non-crossings of Brownian pillow}\\

        \centerline{Enkelejd Hashorva        \footnote{Department of Mathematical  Statistics and Actuarial Science, University of Bern, Sidlerstrasse 5,
        CH-3012 Bern, Switzerland;enkelejd.hashorva@stat.unibe.ch \\
        \today{}}       }
        \centerline{\textsl{University of Bern}}


\end{center}
{\bf Abstract.} Let $B_0(s,t)$ be a Brownian pillow with continuous
sample paths, and let $h,u:[0,1]^2\to\R $ be two measurable functions. In this paper we
derive upper and lower bounds for the boundary non-crossing
probability $\psi(u;h):=\pk{B_0(s,t)+h(s,t) \le u(s,t), \forall s,t\in [0,1]}$. Further we investigate the asymptotic
behaviour of $\psi(u;\gamma h)$ with $\gamma$ tending to $\IF,$ and solve a related minimisation problem.


\section{Introduction}
Let $B_0(s,t), s,t\in [0,1]$ be a Brownian pillow with continuous sample
paths. Its covariance function $K$ is a product of two covariance functions defined by
$$ K((s_1,t_1),(s_2,t_2))= K_1(s_1,t_1)K_2(s_2,t_2), \quad s_i,t_i\in [0,1],i=1,2,$$
with $K_i(s,t)=\min(s,t)-ts , \ i=1,2$ the covariance function of a
Brownian bridge.

Our concern in this article is the boundary non-crossing probability
\BQN \label{eq:def:bnc}
\psi(u;h):=\pk{B_0(s,t)+ h(s,t) \le u(s,t), \forall s,t\in [0,1]},
\EQN
with $h$ a trend function and $u$ a measurable boundary function.\\
When considering Brownian bridge and Brownian motion the corresponding non-crossing probability can be explicitly
calculated if $h$ and $u$ are polygonal lines, see e.g., Goovaerts and Teunen (1993), Wang and P\"otzelberger (1997), Novikov et al.\ (1999),
Janssen and Kunz (2004), Borovkov and Novikov (2005) and the references therein.
Such explicit formulae  are not available in our setup of the multi-parameter processes.

Our novel results presented below are:\\
$a)$ upper and lower bounds for $\psi(u;h)$, $b$) a large deviation type result for the boundary non-crossing probability
$\psi(u;\gamma h)$ with  $\gamma \to \IF$, and $c)$
 we solve a related minimisation problem.

We comment briefly the result mentioned in $b)$.
Given a function $g:[0,\IF)^2\to \R$ we denote by $g''$ its partial derivative  obtained by differentiating both components, provided it exists.
From the large deviation theory (see e.g., Lifshits (1995) or Ledoux (1996))
for any positive constant  $c$ and any trend function  $h:[0,1]^2\to \R$ with a square integrable
partial derivative $h''$ (i.e., $\int_{[0,1]^2} (h''(s,t))^2\, ds dt< \IF$) we obtain
\BQN\label{eq:0}
\lim_{\gamma \to \infty} 2\gamma^{-2} \ln \pb{\sup_{s,t\in [0,1]}(B_0(s,t)+\gamma h(s,t)) \le  c} &= & - \int_{[0,1]^2} (\uhh''(s,t))^2\,
ds dt\in (-\IF,0],
\EQN
with $\uhh$ the solution of the minimisation problem
\BQN\label{min:var:0}
\inf_{g\ge h} \int_{[0,1]^2} (g''(s,t))^2\, ds dt,
\EQN
where the functions $g:[0,1]^2\to \R$ in the minimisation problem are assumed to possess a square integrable
partial derivative $ g''$, and  $g,h$ vanish on the boundary of $[0,1]^2$.

Compared to \eqref{eq:0} our new result is a sharper asymptotic
estimate of the boundary non-crossing probability of interest. In
the special case $h$ being a product of two concave functions
$h_1,h_2:[0,1]\to [0,\IF)$ with $h_i(0)=h_i(1)=0,i=1,2$ we show (see below \eqref{eq:fi})
\BQN\label{eq:show}
\lefteqn{\pb{\sup_{s,t\in [0,1]}(B_0(s,t)+\gamma h_1(s)h_2(t))\le c}}\notag\\
 &=&  \exp\Bigl(-  \frac{\gamma^2}{2}\prod_{i=1,2} \int_{[0,1]}( h_i'(x))^2\, \lambda( dx)
+ c\gamma \prod_{i=1,2} [ h_i'(1)- h_i'(0)]+ z(\gamma)\Bigr),
\EQN
where
$$ -A \gamma^{2/3} \ln ^3\gamma \le  z(\gamma)\le  \ln \pb{\sup_{s,t\in [0,1]}B_0(s,t)\le c}$$
holds for all large $\gamma$ with a positive constant $A$ not depending on $\gamma$. Here $h'_i$
is a right continuous version of the derivative of $h_i,i=1,2$ and $\lambda$ is the
Lebesgue measure on  $[0,1]$.

We derive \eqref{eq:show}  utilising a known small ball result for Brownian pillow. Indeed the small ball problem for both the Brownian pillow and the Brownian sheet is  investigated by several authors, see Kuelbs and Li (1992), Talagrand (1994), Cs\'aki et al.\ (2000), Khoshnevisan and Pemantle (2000), Li and Shao (2001), Koning  and Protasov (2003), Fill and Torcaso (2004), Gao et al.\ (2004), Gao  and Li (2006, 2007) or  Karol' et al.\ (2008) among many other references.

A consequence of the Gaussian shift inequality (see Li and Kuelbs (1998))
and \eqref{eq:show} is the following bound (set $D$ for the set of all concave functions $f:[0,1]\to [0,\IF)$)
\BQN\label{eq:some}
\pb{\sup_{s,t\in [0,1]}B_0(s,t)\le  c} \le \inf_{h\in D}
\Phi\biggl( c^2\biggl(\frac{ h'(1)- h'(0)}{\int_{0}^1 (h'(x))^2\, \lambda(dx)}\biggr)^2\biggr),
\EQN
with $\Phi$ the distribution function of a Gaussian random variable with mean 0 and variance 1.
Since the upper bound in \eqref{eq:some} is not smaller than 1/2, the above inequality is
of some interest provided that  $\psi(0;c)\in (1/2,1)$.

Organisation of the paper: In the next section we present some notation and preliminary results.
The main results are discussed in Section 3. Section 4 explains the simple situation where the trend function $h$ is a product of two trend functions. Proofs of all the results are relegated to Section 5 followed by
a  short Appendix with two results on the Riemann-Stieltjes integral.

\section{Preliminaries}
We introduce first a Hilbert spaces related to the covariance function of the Brownian pillow, which can also be
seen as tensor product of Hilbert spaces related to the covariance function of the Brownian bridge.
Then we provide a result utilised in solving the minimisation problem \eqref{min:var:0}.

The Reproducing Kernel Hilbert Space (RKHS) related to the covariance function of a  Brownian pillow, denoted by $\kHy$, is given by
\BQNY \kHy:=\Biggl\{h:[0,1]^2\to\R~|~\exists
h''\in L_2([0,1]^2,\lambda^2), \mbox{ with } h(s,t)=\int_{[0,s] \times [0,t]}h''(x,y)~\lambda^2(dx,dy),\\
\quad h(0,s)=h(1,s)=h(t,0)=h(t,1)=0, \quad \forall s,t\in [0,1] \Biggr\},
\EQNY
where $L_2([0,1]^2,\lambda^2)$ is the set of all real functions defined on
$[0,1]^2$ which are square integrable with respect to the Lebesgue measure $\lambda^2$ on  $[0,1]^2$. The inner product is
\BQNY
\sprod{h_1}{h_2} &= &\int_{[0,1]^2} h''_1(x,y)h''_2(x,y)\, \lambda^2(dx,dy), \quad
h_1,h_2\in \kHy
\EQNY
and the  corresponding norm of $h\in \kHy$ is $\norm{h}:=\sprod{h}{h}^{1/2}$.\\

As shown in Koning  and Protasov (2003) another approach to deal with $\kHy$ is to construct this Hilbert space  as
the tensor product of two RKHS, i.e.,
 $\kHy= \kHya \otimes \kHya$ with $\kHya$ the RKHS of the covariance function of the Brownian bridge defined by
\BQNY \kHya&:=&\Biggl\{h:[0,1]\to\R~|~\exists
h'\in L_2([0,1],\lambda ), \mbox{ with } h(s)=\int_{[0,s]}h'(x)\, \lambda(dx), \quad h(0)=h(1)=0\Biggr\},
\EQNY
where $L_2([0,1],\lambda )$ is the set of all real functions defined on $[0,1]$ that are square integrable with respect to $\lambda$.
The inner product of $\kHya$ is
$$\sprod{h_1}{h_2} = \int_{[0,1]} h'_1(x)h'_2(x)\, \lambda(dx), \quad h_1,h_2\in \kHya  $$
and the corresponding norm is denoted again by $\norm{\cdot}$.
Any element $h\in \kHy$ can be identified by $h_1,h_2 \in \kHya$ so that $h=h_1 \otimes h_2$ (see Koning  and Protasov (2003)).

In the following for any trend function $h\in \kHy$ we write by $h''$ its right continuous derivative.\\

Lemma 2 in Janssen  and H\"ulya (2008) is crucial for our next result. Define the closed convex sets
$$\CON:=\{h\in \kHy: h(s,t)\le 0, \forall s,t\in [0,1]\}, \quad
\CONB:=\{h\in \kHy: h(s,t)\ge 0, \forall s,t\in [0,1]\},$$
and let $\PCON,\PCONB$ be the polar cones of $\CON$  and $\CONB$, respectively defined by
$$\PCON:= \{h\in \kHy:  \sprod{h}{v} \le 0, \forall v\in \CON\}, \quad
\PCONB:= \{h\in \kHy:  \sprod{h}{v} \ge 0, \forall v\in \CONB\}.$$
Further denote by $BV_H(T), T\subset   \R^2$ the class of functions $f:T \to \R$
which have bounded variation in the sense of Hardy (see e.g., Adams and Clarkson (1934), M\'{o}ricz (2002)).

\BL\label{lem:0}
Let $h\in \kHy$ be a given function and let $\PRJ,\POLPRJ$ be the unique projections of $h$
into  $\CON$ and  the polar cone $\PCON$, respectively. \\
a) If $\POLPRJ''$ is a right continuous partial derivative of $\POLPRJ$ such that
$\POLPRJ'' \in BV_H([0,1]^2)$, then for any function $g:[0,1]^2\to [0,\IF)$
Riemann-Stieltjes integrable with respect to $\POLPRJ''$
the Riemann-Stieltjes  integral $I(g):=\int_{[0,1]^2} g(s,t)\, d \POLPRJ''(s,t) $ satisfies $I(g)\ge 0$.\\
b) We have
\BQN\label{eq:lem:0:1}
h&=&\PRJ+\POLPRJ, \quad \sprod{\PRJ}{\POLPRJ}=0.
\EQN
c) If $h=h_1+h_2$ with $h_1\in \CON, h_2\in \PCON$ such that $\sprod{h_1}{h_2}=0$, then
$h_1=\PRJ$ and $h_2=\POLPRJ$.\\
d) The unique solution $\uhh$ of the minimisation problem
\BQN\label{min:var}
\min_{g\ge h, g\in \kHy} \norm{g}
\EQN
is $\uhh=\POLPRJ$ satisfying  further $\norm{\uhh}= \min \{\norm{g}: g\in \PCON, g\ge h \}$.
\EL
\def\uhhV{\uhh}

We note in passing that a similar decomposition to \eqref{eq:lem:0:1} can be stated for $h\in \kHy$
in terms of the unique projections  $\PRJB, \POLPRJB$ of $h$ into  $\CONB$ and  the polar cone $\PCONB$, respectively.
Furthermore, $b),c)$ hold for some general Hilbert space.\\
We write alternatively $\uhh, \lhh$ instead of $\POLPRJ,\POLPRJB$. The above lemma immediately implies
\BQN \label{nromhGEnormUHH}
 \lhh(s,t) \le h(s,t) \le \uhh(s,t), \quad  \forall s,t\in [0,1], \quad \text{ and }
\norm{h} &\ge& \max( \norm{\uhh}, \norm{\lhh}), \quad \forall h\in \kHy.
\EQN
Furthermore, for any two functions $h,q\in \kHy, q \ge h$
 \eqref{min:var:0}  and Lemma \ref{lem:0} yield
\BQN\label{eq:order:norms}
\norm{q} &\ge& \norm{h},
\EQN
provided that $\uhh=h, \underline{q}=q$.

\section{Main Results}
Let $B_0(s,t),s,t\in [0,1]$ be a Brownian pillow with continuous sample paths, and let $h\in \kHy$ be
a given trend function. For some measurable boundary function $u:[0,1]^2\to \R$ we define  the
boundary non-crossing probability $\psi(u;h)$ as in \eqref{eq:def:bnc}.
Throughout the rest of the paper we assume that $\psi(u;0)\in (0,1)$. Since $h\in \kHy$ the Cameron-Martin formula
(see e.g., Kuelbs (1976), Lifshits (1995), Li and Kuelbs (1998) or Li and Shao (2001)) implies
\BQN\label{eq:CM}
\psi(u;h)&=& \exp\Bigl(-  \frac{1}{2}\norm{h}^2\Bigr)
\Eg{\exp\Bigl(\int_{[0,1]^2} h''(s,t)\,d B_0(s,t)  \Bigr)\vk{1}(B_0(s,t) \le u(s,t), \forall s,t \in [0,1])},
\EQN
where $\vk{1}(\cdot)$ is the indicator function.

Li and Kuelbs (1998) show that the Cameron-Martin translation implies important shift inequalities for some general Gaussian processes.
Applying their Theorem 1' we have
\BQN\label{eq:WL}
\Phi(\theta - \norm{h}) \le \psi(u;h) \le \Phi(\theta + \norm{h}),
\EQN
where $\Phi$ is the Gaussian distribution function on $\R$ with mean 0 and variance 1, and $\theta$ is such that
$\Phi(\theta)= \psi(u;0)$. When $\norm{h}$ is small the lower and the upper bounds in \eqref{eq:WL}
are close to the non-crossing probability of interest, since
$ \lim_{\gamma \to 0}\psi(u;\gamma h)= \psi(u;0)= \Phi(\theta)$.  As  $\gamma \to \IF$ the upper bound in \eqref{eq:WL} tends to 1,
 whereas the lower bound and $\psi(u; \gamma h)$ tend to 0. Note in passing that as in P\"otzelberger and Wang (2001)  we obtain
\BQN\label{eq:potz:difh0}
\abs{\psi(u;\gamma h)- \psi(u;0)} &\le & 2\Phi(\gamma \norm{h}/2) -1 \le \frac{\gamma \norm{h}}{\sqrt{2 \pi}}, \quad \forall \gamma \in (0,\IF).
\EQN
One important criteria which we will look at when discussing bounds for the non-crossing probability of interest is
their performance for both small or large trend functions. In our first result below we provide upper and lower bounds for the boundary non-crossing
probability $\psi(u;h)$. If we consider further the trend function $\gamma h$,
then the bounds perform well when $\gamma \to 0$.
\begin{theorem}\label{theo:00}
Let $h,u:[0,1]^2 \to \R$ be two measurable functions such that $\psi(u;0)\in (0,1)$. If $ h\in \kHy$, then we have
\BQN\label{eq:00:1}
\Phi(\theta- \norm{\uhh}) \le &\psi(u;h) \le& \Phi(\theta+ \norm{\lhh}), \quad \theta:= \Phi^{-1}(\psi(u;0)),
\EQN
with $\uhh$, $\lhh$ as defined in Section 2 and $\Phi^{-1}$ the inverse of $\Phi$. Furthermore
\BQN\label{eq:00:2}
- \frac{\norm{\uhh}}{\sqrt{2 \pi}} \le & \psi(u;h) - \psi(u;0) &\le \frac{\norm{\lhh}}{\sqrt{2 \pi}}.
\EQN
\end{theorem}
When $h\not= \uhh$ or $h \not= \lhh$  in view of \eqref{nromhGEnormUHH} we see that
\eqref{eq:00:2} yields better bounds than \eqref{eq:potz:difh0}.
 By \eqref{eq:00:2} we obtain
\BQN\label{eq:00:2b}
- \gamma \frac{\norm{\uhh}}{\sqrt{2 \pi}} \le & \psi(u;\gamma h) - \psi(u;0) &\le \gamma \frac{\norm{\lhh}}{\sqrt{2 \pi}} , \quad
\forall \gamma>0,
\EQN
which is of some interest  when $\gamma$ tends to 0 since both the lower and the upper bounds converge to 0.

As mentioned in the Introduction if $\gamma$ tends to infinity, then we have the logarithmic asymptotic behaviour
\BQN\label{eq:large:2}
\lim_{\gamma \to \infty} 2\gamma^{-2} \ln \psi(u;\gamma h)&= & - \norm{\uhhV}^2, \quad \forall h\in \kHy,
\EQN
with $\uhhV$ the unique solution of the minimisation problem $\eqref{min:var}$.\\
Next, we derive explicit upper and lower bounds for $\psi(u;h)$, which perform better asymptotically (for trend
function becoming large) as those implied by \eqref{eq:00:1}.

\begin{theorem}  \label{th1}
Let  $h\in \kHy$ be a given trend function, and let  $u, l:[0,1]^2\to \R$ be two measurable functions.
If the partial derivative $\uhh''$ of the projection of $h$ into
its polar cone satisfies $\uhh''\in \BVHA$ and is right continuous, then
\BQN \label{last:1}
\uhhV&:= &\inf_{g\ge h, g\in \PCON, g\in \BVHA} g
\EQN
and further $\uhh$ is the smallest majorant
of $h$ such that its right continuous partial derivative belongs to $\BVHA$ and generates a finite positive measure.\\
Moreover, if the Riemann-Stieltjes integral
$\int_{[0,1]^2} v(s,t)\,  d \uhh''(s,t)$ is finite for both  $v=l$ and $v=u$ and $\psi(u;0)\in (0,1)$, then
\begin{eqnarray}\label{eqa:main:1}
\psi(u;h)& \le & \psi(u;\deltahH)\exp \Bigl( -\frac{1}{2}
\norm{\uhh}^2+  \int_{[0,1]^2} u(s,t)\,  d \uhh''(s,t) \Bigr)
\end{eqnarray}
and
\begin{eqnarray}\label{eqa:main:1b}
\psi(u;h)& \ge & \pk{ l(s,t) \le B_0(s,t)\le u(s,t), \forall s,t\in[0,1]}\notag \\
 &&\times\exp \Bigl( -\frac{1}{2}
\norm{\uhh}^2+  \int_{[0,1]^2} l(s,t)\,  d \uhh''(s,t) \Bigr)
\end{eqnarray}
are valid.
\end{theorem}
{\bf Remarks:} {\it
a) If  $u(s,t):=c\in (0,\IF),\forall s,t\in[0,1],$ then \eqref{eqa:main:1} implies
  \begin{eqnarray} \label{exa:main}
\psi(c;h)& \le & \psi(c;\deltahH )\exp \Bigl( -\frac{1}{2}
\norm{\uhh}^2+  c[\uhh''(1,1) - \uhh''(1,0)-\uhh''(0,1)+ \uhh''(0,0)] \Bigr).
\end{eqnarray}
A lower bound for $\psi(c;h)$ is derived utilising \eqref{eqa:main:1b} with $l(s,t):=-c,\forall s,t\in [0,1]$.\\

b) As in the proof of Proposition \ref{th1} it can be shown that if the trend function
$h\in \kHy$ is such that its right continuous partial derivative $h''$ satisfies $h''\in \BVHA$
and furthermore $h''$ generates a positive measure on $[0,1]^2$,
then the unique solution of the minimisation problem
\eqref{min:var} is $\uhhV=h$. \\

c) An upper bound for $\psi(u;h)$ is the discrete boundary non-crossing probability
$$ \psi_n(u;h):= \pk{ B_0(s_i,t_i)+ h(s_i,t_i) \le u(s_i,t_i), \forall (s_i,t_i)\in T_n}, $$
with $T_n:=\{ (s_i,t_i),i=1, \dots , n \} \subset [0,1]^2$. Hashorva (2005b) shows the asymptotic behaviour
(considering the Brownian bridge) of the corresponding discrete boundary non-crossing probability}.\\

Next, we discuss the asymptotic behaviour of  $\psi(u;\gamma h)$ when $\gamma \to \IF$. Exact asymptotics of the non-crossing probabilities of
 the Brownian motion with trend is derived in Hashorva (2005a) which was motivated by
 a large deviation type result obtained in Bischoff et al.\ (2003).
 As in Bischoff et al.\ (2005) we expect that our novel asymptotic result
will have some implications for statistical applications.

\begin{theorem}  \label{th2}
Let  $h,\uhh, u$ be as in Proposition \ref{th1}. Suppose that there exist functions
$u_\varepsilon \in \kHy,\varepsilon >0 $ such that $\norm{u_\varepsilon}=O(1/\varepsilon)$ and
\BQN
\lim_{\varepsilon \to 0} u_\varepsilon(s,t)= u(s,t), \quad  u_\varepsilon (s,t) \le  u(s,t)- \varepsilon, \quad \forall s,t\in [0,1].
\EQN
If the Riemann-Stieltjes integral $I_\epsilon:=\int_{[0,1]^2}u_\epsilon(s,t)\,  d \uhh''(s,t)$
exists  and $\abs{I_\epsilon} \le M \in (0,\IF), \forall \epsilon>0$, then
\BQN
\lim_{\epsilon \to 0} I_\epsilon &=&I:= \int_{[0,1]^2} u(s,t)\,  d \uhh''(s,t),
\quad \abs{I}\le M,
\EQN
and
\BQN\label{eqa:main2:2}
\psi(u;\gamma h) &=&
\exp\Bigl(-  \frac{\gamma^2}{2}\norm{\uhh}^2 + \gamma \int_{[0,1]^2} u(s,t)\, d\uhh''(s,t)
+ z(\gamma)\Bigr)
\EQN
holds, where for all large $\gamma$
\BQN\label{eqAa}
 -A \gamma^{2/3} \ln ^3\gamma \le  z(\gamma)\le  \ln  \pk{B_0(s,t)\le u(s,t),
\forall s,t\in [0,1]: \uhh(s,t) = h(s,t)},
\EQN
with  a positive constant $A$ not depending on $\gamma$.
\end{theorem}
In view of the above asymptotics and \eqref{eq:00:1} we obtain a simple upper bound for $\psi(u;0)$.
\begin{corollary}
Let $u:[0,1]^2\to \R$ be a measurable function satisfying the assumptions of Proposition \ref{th2}.
Then we have
\BQN \label{eq:u:o:korr}
\psi(u;0) & \le & \inf_{h\in \kHy, h'' \in \BVHA: \norm{\uhh}>0}\Phi\biggl( \norm{\uhh}^{-1}\int_{[0,1]^2} u(s,t)\, d\uhh''(s,t)\biggr).
\EQN
\end{corollary}

\def\ve{\varepsilon}

\bigskip

{\bf Remarks:}
{\it
a) If the function $u$ in Proposition \ref{th2} satisfies  $u(s,t)> \mu\in (0,\IF),
\forall s,t \in [0,1]$ where  $(s,t)$ belongs to the boundary of $[0,1]^2$, and there exist functions $w_\ve: [0,1]^2 \to \R, \ve>0$ such that
$uw_\ve\in \kHy,\ve >0$, then we may define $u_\ve$ in Proposition \ref{th2} by $u_\ve:= u w_\ve - \epsilon,\epsilon>0 $.
When $u$ is a positive constant, then functions $u_\ve,\ve >0$ satisfying the assumption of Proposition \ref{th2} can be easily constructed.
If  $u_\ve$ is continuous, then the Riemann-Stieltjes integral $I_\ve:= \int_{[0,1]^2} u_\ve(s,t)\, d \uhh'' (s,t)$ in Proposition \ref{th2} is finite.

b) When  $\uhh''$ is almost surely continuous with respect to the Lebesgue measure $\lambda^2$, then
instead of assuming that $\uhh$ has a bounded variation in the sense of Hardy (Lemma \ref{lem:0}, Proposition \ref{th1} and \ref{th2})
we may impose the weaker assumption $\uhh$  has a bounded variation in the sense of Vitali (see Appendix below and Lemma \ref{lem:mo:2}).

c) Our results can be easily extended to the $d$-dimensional setup considering $B_0(s_1, \dots, s_d), s_i \in [0,1], i\le d$
a Brownian pillow with continuous sample paths. The term $\ln^3 \gamma$ in \eqref{eqAa} should then be replaced by  $\ln ^{2d-1} \gamma$.

d) Similar results can be stated for considering instead of $B_0$ a Brownian sheet $B(s,t), s,t \in [0,\IF)$ with continuous sample paths.
For instance  Proposition \ref{th1} holds with $\underline{h}$ the solution of the minimisation problem  (1.3) where $g,h$ have square integrable partial derivatives satisfying further  $g(0,s)=h(0,s)=h(t,0)=g(t,0)=0,s,t\in [0,\IF)$.
}

\section{Product Trend Functions}
As demonstrated in the previous section the non-crossing probability $\psi(u;h)$
can be bounded by some functions which depend on the solution of the minimisation problem \eqref{min:var}. We discuss below
an instance where the solution of \eqref{min:var} can be easily determined. Let therefore  $h_1,h_2\in \kHya$, and let $B_0(s),s\in [0,1]$ denote a Brownian bridge with continuous sample paths. If $u_1,u_2:[0,1] \to \R$ are two measurable functions with $u_i(0),u_i(1)>0,i=1,2$, then  we have (see Bischoff and Hashorva (2005))
\BQNY
\lefteqn{\pk{ B_0(s)+h_i(s) \le u_i(s), \forall s\in [0,1]}}\\
&\le&
\pk{ B_0(s) \le  u_i(s)+ \widetilde h_i(s)- h_i(s), \forall s\in [0,1]}\exp\Bigl(-  \frac{1}{2}\norm{\widetilde h_i}^2
+ \int_{[0,1]} u_i(s) \,  d (- \widetilde h'_i(s))\Bigr),
\EQNY
where $\widetilde h_i,i=1,2$ is the smallest concave majorant of $h_i$ and $\widetilde h'_i$
is a right continuous derivative of $\widetilde h_i$. Furthermore, $\widetilde h_i$ is the unique solution of the minimisation problem
\BQN \label{min:var:3}
\min_{g\in \kHya, g\ge h_i}\norm{g}, \quad i=1,2.
\EQN
Set in the following $h(s,t):=h_1(s)h_2(t), \tilde h(s,t):= \tilde h_1(s)\tilde h_2(t), s,t\in [0,1]$, and write $h= h_1 \times h_2, \tilde h= \tilde h_1 \times \tilde h_2$. In the next lemma we show that for special trend functions the unique solution of \eqref{min:var} with
$h=h_1\times h_2 \in \kHy$ is simply $\widetilde  h$.

\begin{lemma}  \label{th6}
Let  $h:= h_1 \times h_2, h_1,h_2 \in \kHya$, and denote by $\widetilde  h_i,i=1,2$ the smallest concave
majorant of $h_i,i=1,2$. If
\BQN\label{eq:unb1}
\widetilde  h (s,t) \ge h(s,t) , \quad \forall s,t \in [0,1],
\EQN
then the unique solution $\uhhV$ of \eqref{min:var:0} is $\uhh:= \widetilde  h$.
\end{lemma}
Clearly, \eqref{eq:unb1} holds if $h_1,h_2$ are both non-negative functions. In the special case that also $u$ is a product function
 we have the following immediate result.

\BK\label{eq.frauenk}
Let  $h_i,\widetilde h_i,i=1,2$  satisfy the assumption of Lemma \ref{th6}, and let
 $u_i,l_i:[0,1]\to \R, i=1,2$ be  measurable functions. If  the Riemann-Stieltjes integral $\int_{[0,1]} v_i(s)\, d (- \widetilde h_i'(s))$
is a finite constant for $i=1,2$ and $v_i=l_i$ or $v_i=u_i$,  then we have
  \begin{eqnarray}
\label{eq:th6:1}
\psi(u ;h)  \le  \psi(u;  h-\tilde h ) \exp \Bigl( -\frac{1}{2}
\norm{\widetilde h_1}^2\norm{\widetilde h_2}^2+  \prod_{i=1,2} \int_{[0,1]} u_i(s)\,  d (- \widetilde  h_i'(s)) \Bigr),
\end{eqnarray}
with $h:= h_1 \times h_2, \tilde h:= \tilde h_1 \times \tilde h_2, u:= u_1 \times u_2$, and further
  \begin{eqnarray}
\psi(u;h)& \ge & \pk{ l_1(s)l_2(t)\le B(s,t)\le u_1(s) u_2(t), \forall s,t \in [0,1]}\notag \\
&& \times  \exp \Bigl( -\frac{1}{2}
\norm{\widetilde h_1}^2\norm{\widetilde h_2}^2
+  \prod_{i=1,2} \int_{[0,1]} l_i(s)\,  d (- \widetilde  h_i'(s)) \Bigr).
\EQN
\EK

\BK\label{kor:fi}
Under the assumptions and the notation of Corollary \ref{eq.frauenk} if further
$\min_{s\in [0,1]}u_i(s)> C\in (0,\IF), i=1,2$ and $u_i, i=1,2$ are absolute continuous with $u'_i$ satisfying
$\int_{[0,1]}(u_i'(s))^2\, \lambda(ds)< \IF$, then we have
\BQN\label{eq:fi}
\psi(u_1\times u_2;\gamma h_1\times h_2) &=&  \exp\Bigl(-  \frac{\gamma^2}{2}
\norm{\widetilde h_1}^2\norm{\widetilde h_2}^2 + \gamma \prod_{i=1}^2\int_{[0,1]} u_i(s)\, d (-\widetilde h_i'(s)
)+ z(\gamma)\Bigr)
\EQN
with $z(\gamma)$ satisfying
$$ -A  \gamma^{2/3} \ln ^3\gamma \le  z(\gamma)\le
\ln  \pk{B_0(s,t) \le  u_1(s)u_2(t), \forall s,t\in [0,1]: \widetilde h_1(s)\widetilde h_2(t)=  h_1(s) h_2(t)}$$
for all large $\gamma$ where $A$ is a positive constant not depending on $\gamma$. Furthermore
  \begin{eqnarray}
\psi(u;0)& \le & \inf_{h_1,h_2 \in \kHya: \norm{\widetilde  h_1}\norm{\widetilde  h_2}>0} \Phi\biggl( (\norm{\widetilde  h_1}\norm{\widetilde  h_2})^{-1}
\prod_{i=1,2} \int_{[0,1]} u_i(s)\,  d (- \widetilde  h_i'(s)) \biggr).
\EQN
\EK


\section{Proofs}
{\bf Proof of Lemma \ref{lem:0}}: Let $g,h \in \kHy$ be two given functions.
If $h'' \in \BVHA$ with $h''$ a right continuous partial derivative of $h$,
then we have by \eqref{eq:int:p2} and the integration by parts formula (see Lemma 2 and Lemma 3 in M\'{o}ricz (2002) and \eqref{eq:int:p})
\BQN \label{eq:biv:partial}
\sprod{ g}{h} &=&\int_{[0,1]^2} g''(s,t)  h''(s,t) \, \lambda^2(ds, dt)\notag \\
&=&\int_{[0,1]^2} h''(s,t) \, d g(s,t)\notag \\
&=&\int_{[0,1]^2} g(s,t) \, d h''(s,t).
\EQN
Consequently, for any $g\in \CON$ by the assumption on $\POLPRJ''$  we have $\sprod{ g}{ \POLPRJ} \le 0.$
Hence for any function $g:[0,1]^2\to [0,\infty)$ which is Riemann-Stieltjes integrable
with respect to $\POLPRJ''$  on $[0,1]^2$
we have for the corresponding Riemann-Stieltjes integral
\BQN
\label{eq:biv:partial:2}\int_{[0,1]^2} g(s,t) \, d \POLPRJ''\ge 0.
\EQN
The proof of statements $b)$ and $c)$  follows  immediately by Lemma 2 in Janssen  and H\"ulya (2008).\\
We show next statement $d)$. Let $\tilde h\in \kHy$ be a given function such that $\tilde h:= g+ h$ with $g(s,t)\ge 0,\forall s,t\in [0,1]$. By the properties of $\POLPRJ$ we have $\sprod{\POLPRJ}{g}\ge 0$, hence we may write
\BQNY
\norm{\tilde h}^2&=& \norm{g+ h}^2\\
&=&\norm{\POLPRJ+ g+h-\POLPRJ}^2\\
&=&\norm{\POLPRJ}^2+ 2\sprod{\POLPRJ}{g+h-\POLPRJ}+ \norm{g+h-\POLPRJ}^2\\
&=&\norm{\POLPRJ}^2+ 2\sprod{\POLPRJ}{g}+ 2\sprod{\PRJ}{\POLPRJ}+\norm{g+h-\POLPRJ}^2\\
&=&\norm{\POLPRJ}^2+ 2\sprod{\POLPRJ}{g}+ \norm{g+h-\POLPRJ}^2\\
&\ge &\norm{\POLPRJ}^2+ 2\sprod{\POLPRJ}{g}\\
&\ge &\norm{\POLPRJ}^2.
\EQNY
Since further $\POLPRJ(s,t)\ge h(s,t), \forall s,t\in [0,1]$ it follows that the solution of the minimisation problem \eqref{min:var} is
$\POLPRJ$. Clearly, its solution is unique, thus the result follows. \hfill \Box\\

{\bf Proof of Proposition \ref{theo:00}:} By \eqref{nromhGEnormUHH} and \eqref{eq:WL} we see that \eqref{eq:00:1} follows easily.
The proof of \eqref{eq:00:2} can be established along the lines of the proof of Lemma 5
 in Janssen  and H\"ulya (2008), thus the result.  \hfill \Box\\

{\bf Proof of Proposition \ref{th1}}: Let $\CON,\PCON$ be as in Section 2, and let
$\POLPRJ$ be the projection of $h$ into the polar cone $\PCON$.
In view of statement $b)$ of Lemma \ref{lem:0}
$$ h= \PRJ+ \POLPRJ, \quad \norm{h}^2= \norm{\POLPRJ}^2+\norm{\PRJ}^2. $$
Furthermore, $\psi(u;h) \ge \psi(u;\POLPRJ)$. Next, applying the Cameron-Martin formula we obtain  (set
$\njTAU:=\vk{1}(B_0(s,t)  \le  u(s,t),\forall s,t\in [0,1])$
\BQNY
\psi(u;h) & = & \exp\biggl(-  \frac{1}{2}\norm{h}^2  \biggr)
\Eg{ \exp\biggl( \int_{[0,1]^2} h''(s,t)\, dB_0(s,t) \biggr)\njTAU} \\
& = & \exp\biggl(-  \frac{1}{2}\norm{\POLPRJ}^2 \biggr)
\Eg{ \exp\biggl(-  \frac{1}{2}\norm{\PRJ}^2 +   \int_{[0,1]^2} \PRJ''(s,t)\, dB_0(s,t)\\
&&
+\int_{[0,1]^2} \POLPRJ''(s,t)\, dB_0(s,t) \biggr)\njTAU}.
\EQNY
Since $\POLPRJ''\in \BVHA$ is right continuous and $B_0(s,t)$ has continuous sample paths by the integration by  parts formula
\eqref{eq:int:p} for the Riemann-Stieltjes integral we have almost surely
$$  \int_{[0,1]^2}  B_0(s,t) d \POLPRJ''(s,t) = \int_{[0,1]^2} \POLPRJ''(s,t)\, dB_0(s,t). $$
Consequently,  we may further write (recall \eqref{eq:biv:partial:2})
\BQNY
\psi(u;h) & = & \Eg{ \exp\biggl(-  \frac{1}{2}\norm{\PRJ}^2 +   \int_{[0,1]^2} \PRJ''(s,t)\, dB_0(s,t)
 +\int_{[0,1]^2}  B_0(s,t) d \POLPRJ''(s,t) \biggr)\njTAU} \\
& \le  & \exp\biggl(-  \frac{1}{2}\norm{\POLPRJ}^2 +\int_{[0,1]^2}  u(s,t) d \POLPRJ''(s,t)\biggr)\\
&&
\Eg{ \exp\biggl(-  \frac{1}{2}\norm{\PRJ}^2 +   \int_{[0,1]^2} \PRJ''(s,t)\, dB_0(s,t)\biggr)\njTAU} \\
& =  & \exp\biggl(-  \frac{1}{2}\norm{\POLPRJ}^2 +\int_{[0,1]^2}  u(s,t) d \POLPRJ''(s,t)\biggr)
\psi(u;\PRJ ).
\EQNY
Clearly, by the definition $\psi(u;h)\ge \psi(u;\POLPRJ)$. Applying \eqref{eq:large:2} to
$\psi(u; \gamma \POLPRJ), \gamma >0$ we find
\BQNY
\ln \psi( u; \gamma h)&=& - (1+o(1))\frac{\gamma ^2}{2} \norm{\POLPRJ}^2, \quad \gamma \to \IF,
\EQNY
hence by \eqref{eq:large:2} the unique solution of \eqref{min:var} equals $\POLPRJ$.
Since $\POLPRJ\ge h$ and $\POLPRJ\in \PCON$, then $\uhh=\POLPRJ$ and \eqref{last:1} follows. \\
We show next the last claim \eqref{eqa:main:1b}. Utilising again the Cameron-Martin formula we have
\BQNY
\psi(u;h) &\ge & \psi(u;\uhh)\\
&\ge &\pk{ l(s,t)\le B_0(s,t)+\uhh(s,t) \le  u(s,t),\forall s,t\in [0,1]}\\
&=& \exp\Bigl(-  \frac{1}{2}\norm{\uhh}^2\Bigr)
\Eg{\exp\Bigl(\int_{[0,1]^2} \uhh''(s,t)\,d B_0(s,t)  \Bigr)
\vk{1}(l(s,t)\le B_0(s,t) \le  u(s,t),\forall s,t\in [0,1])}\\
&=& \pk{ l(s,t)\le B_0(s,t)\le u(s,t),\forall s,t\in [0,1]}\exp\Bigl(-  \frac{1}{2}\norm{\uhh}^2+\int_{[0,1]^2} l(s,t)\, d \uhh''(s,t)\Bigr),
\EQNY
hence the proof is established. \hfill \Box\\

{\bf Proof of Proposition \ref{th2}}: Set next
$$ \we(s,t):=\uhh(s,t)-  u_\epsilon(s,t), \quad \forall s,t \in [0,1].$$
Applying Cameron-Martin formula we obtain
\BQNY
\psi(u;h) &\ge & \psi(u;\uhh) \\
&= &\pk{ B_0(s,t)+\uhh(s,t) \le u(s,t),\forall s,t\in [0,1]}\\
&\ge &\pk{ B_0(s,t)+\uhh(s,t) \le u_\epsilon(s,t) +\epsilon ,\forall s,t\in [0,1]}\\
&> & \exp\Bigl(-  \frac{1}{2}\norm{\we}^2\Bigr)
\Eg{\exp\Bigl(\int_{[0,1]^2} \we''(s,t)\,d B_0(s,t)  \Bigr) \vk{1}(-\epsilon \le B_0(s,t) \le \epsilon ,\forall
s,t\in [0,1])}.
\EQNY
Define the Gaussian random variable $$ Z: = \int_{[0,1]^2} \we''(s,t)\,d B_0(s,t).$$ 
Clearly, $Z$ has mean 0 and variance $\norm{\we}^2$. For $\varepsilon>0$ small enough we
have $\norm{\we}\in (0,\IF)$. For any constant $C\in \R$ and $\varepsilon$ small enough we may write
\begin{eqnarray*}
\lefteqn{\E{\exp( Z) \nj(    -\epsilon \le  B_0(s,t) \le \epsilon, \forall s,t \in [0,1])}}\\
&=&\E{\exp( Z) \nj(    -\epsilon \le  B_0(s,t) \le \epsilon, \forall s,t \in [0,1])[\nj(Z< C)+ \nj(Z\ge C)]}\\
&\ge&\E{\exp( Z) \nj(    -\epsilon \le  B_0(s,t) \le \epsilon, \forall s,t \in [0,1]) \nj(Z\ge C)}\\
&\ge&\exp(C)\pk{     -\epsilon \le  B_0(s,t) \le \epsilon, \forall s,t \in [0,1], Z\ge C}\\
&=&\exp(C)\Bigl[\pk{ \sup_{s,t\in [0,1]}\abs{B_0(s,t)}< \epsilon}-
\pk{   -\epsilon \le  B_0(s,t) \le \epsilon, \forall s,t \in [0,1],  Z< C}\Bigr]\\
&\ge&\exp(C)\Bigl[ \pk{ \sup_{s,t\in [0,1]}\abs{B_0(s,t)}< \epsilon} - \pk{Z\le C}\Bigr]\\
&=&\exp(C)[ \pk{ \sup_{s,t\in [0,1]}\abs{B_0(s,t)}< \epsilon} - \Phi(C/\norm{\we})].
\end{eqnarray*}
By the small ball asymptotic result (see Fill and Torcaso (2004), Gao and Li (2006, 2007),
Karol' et al.\ (2008)) we have
$$ \pb{ \sup_{s,t\in [0,1]}\abs{B_0(s,t)}< \epsilon}  \ge \exp\biggl(- K \frac{\ln^3(1/\epsilon)}{\epsilon^2}\biggr)$$
for some positive constant $K$ and all $\epsilon> 0$ small enough.
Since
$$\norm{\we}^2= \norm{\uhh}^2- 2 \int_{[0,1]^2} u_\varepsilon(s,t) \, d \uhh''(s,t)+ \norm{u_\varepsilon}^2=O(1/\varepsilon^2)
$$
choosing $C:=-K_* \norm{\we}\ln^{3/2}(1/\epsilon)/\epsilon, K_*\in (0,\IF), K_*^2>K$ and using the Mills-Ratio asymptotics for Gaussian random variables for all $\epsilon >0$ small enough and some positive constants  $c_1,c_2$ we have
\begin{eqnarray*}
\E{\exp( Z) \nj(    -\epsilon \le  B_0(s,t) \le \epsilon, \forall s,t \in [0,1])}&\ge&
\exp\biggl(- \frac{c_1}{\epsilon} - \frac{ c_2 \ln^3(1/\epsilon)}{\epsilon^2}\biggr),
\end{eqnarray*}
implying thus
\BQNY
\psi(u;h) &\ge &  \exp\Bigl(-  \frac{1}{2}\norm{\uhh}^2 + \int_{[0,1]^2} u_\epsilon(s,t)
\, d\uhh''(s,t)  - \frac{c_1}{\epsilon} - \frac{c_2 \ln^3(1/\epsilon)}{\epsilon^2}\Bigr).
\EQNY
Recalling that $\lim_{\varepsilon \to 0} u_\epsilon(s,t)= u(s,t),\forall s,t\in [0,1] $ and $\norm{u_\epsilon}^2=O(1/\varepsilon^2)$
we obtain utilising the result of Proposition \ref{th1} (set next  $\epsilon := \gamma^{-1/3}, \gamma >0$)
\BQNY
\psi(u;\gamma h) &=&  \exp\Bigl(-  \frac{\gamma^2}{2}\norm{\uhh}^2 + \gamma
I+ z(\gamma)\Bigr), \quad \gamma \to \IF,
\EQNY
where $\abs{I}\le M$ with $I:= \int_{[0,1]^2} u(s,t)\, d\uhh''(s,t)$ and
$$ -A \gamma^{2/3} \ln ^3 \gamma \le  z(\gamma)\le  \ln \pb{B_0(s,t)\le u(s,t), \forall s,t\in [0,1]:
 \uhh (s,t)= h(s,t) } $$
is satisfied for all $\gamma$ large and  $A$ a positive constant not depending on $\gamma$.
Hence the result follows. \hfill \Box\\

{\bf Proof of Lemma \ref{th6}:} Set $V:=\{ h\in \kHy: h(s,t) \le 0, \forall s,t\in [0,1]\}$ and
$\uhh:=\widetilde h_1 \times \widetilde h_2$. By the assumptions the function
$g:= \uhh - h_1\times  h_2$ belongs to $V$. Furthermore, for any $v\in V$ we have
$$ \sprod{ v}{ \uhh} = \int_{[0,1]^2} v(s,t) \, d( \widetilde h_1'(s) \widetilde h_2'(t)) \le 0.$$
Consequently $\uhh$ belongs to the polar cone $\PCON$ of $V$.  In view of statement $c)$ in Lemma \ref{lem:0}
 the proof follows if we show that $g$ is orthogonal to  $\uhh$. Since $\widetilde h_i- h_i$ is orthogonal to  $\widetilde h_i,i=1,2$
(see Bischoff and Hashorva (2005)) we have
\BQNY
\sprod{ g}{ \uhh} &=& \sprod{\widetilde h_1 \times \widetilde h_2- h_1 \times h_2}
{ \widetilde h_1 \times \widetilde h_2 }\\
&=& \sprod{\widetilde h_1 \times (\widetilde h_2-h_2)}{ \widetilde h_1 \times \widetilde h_2} -
\sprod{ (\widetilde h_1 - h_1)\times  h_2}{ \widetilde h_1 \times \widetilde h_2 }\\
&=&0,
\EQNY
hence the result follows. \hfill \Box\\

{\bf Proof of Corollary \ref{kor:fi}}: The proof follows easily by the assumptions on $u_i,i=1,2$. \hfill \Box\\

\section{Appendix}
In this short section we provide two results for the
Riemann-Stieltjes integral. \\
Let $f:[0,1]^2 \to \R$ be a given function. If
$f(s,t)=g(s,t)+g_1(s)+g_2(t)$ with $g\in  \BVHA$  and $g_1,g_2$ two
other functions, then $h$ has bounded variation in the sense of
Vitali (write $f\in \BVVA$). In fact $f$ can be expressed as the
difference of two real functions defined on  $[0,1]^2$ which generate a positive
measure on $[0,1]^2$. Thus the class of functions with bounded
variation in the sense of Vitali consists of all real functions defined on $[0,1]^2$ generating a
finite signed  measure.

If $g:[0,1]^2\to \R$ is continuous, then it is well-known that the
Riemann-Stieltjes integral $\int_{[0,1]^2} g(x,y)\, d f(x,y)$
exists, provided that $f\in \BVVA$.  In the next lemma we present  an
integration by parts formula, the case $f \in \BVHA$ is discussed  in
Lemma 1 in M\'{o}ricz (2002).

\begin{lemma}\label{lem:mo:1} Let $f,g:[0,1]^2\to \R$ be two given
function. If $g$ is continuous such that $g(s,t)=0$ for all $(s,t)$
in the boundary of $[0,1]^2$ and $f\in \BVVA$, then the integration by parts formula for the Riemann-Stieltjes integral reads
\BQN\label{eq:int:p}
\int_{[0,1]^2} g(x,y) \, d f(x,y) &=&
\int_{[0,1]^2} f(x,y) \, d g(x,y).
 \EQN
\end{lemma}
{\it Proof:} The proof follows with similar arguments as in Lemma 2 in
M\'{o}ricz (2002), since the four single sums in
the expression (3.8) therein  are equal to 0 due to the fact that $g$
vanishes on the boundary of $[0,1]^2$. \hfill \Box\\

\begin{lemma}\label{lem:mo:2} Let $f,g:[0,1]^2\to \R$ be two given
functions. Assume that $g$ is absolute continuous with
$g(s,t)=\int_{[0,s]\times [0,t]} h(x,y)\, \lambda^2(dx, dy),s,t\in [0,1]$.
If $f\in \BVVA$ and $f$ is almost surely continuous with respect to $\lambda^2$, then we have
\BQN\label{eq:int:p2}
\int_{[0,1]^2} g(x,y) \, d f(x,y) &=&
\int_{[0,1]^2} f(x,y) h(x,y)\, d \lambda^2(dx,dy).
 \EQN
\end{lemma}
{\it Proof:} The proof follows with similar arguments as in Lemma 3 in M\'{o}ricz (2002).  \hfill \Box\\

 {\bf Acknowledgement:} I would like to thank a Referee and Professor Wembo Li for several corrections and suggestions,
Professor M\'{o}ricz for sending \cite{M2002},  Professors Muhammad Aslam Noor and Wolfgang Bischoff for some insights on Hilbert spaces.

\bibliographystyle{plain}

\end{document}